\documentclass[a4paper,UKenglish,cleveref,autoref,thm-restate]{lipicspreprint-v2021}

\title{Mixed-Integer Optimization for Loopless Flux Distributions in Metabolic Networks} 

\titlerunning{Mixed-Integer Optimization for Loopless Fluxes} 

\author{Hannah Troppens}{Zuse Institute Berlin, Germany}{hannahtroppens@gmail.com}{https://orcid.org/0000-0002-8079-9603}{}
\author{Mathieu Besançon}{Université Grenoble Alpes, Inria, LIG, CNRS, France}{mathieu.besancon@inria.fr}{https://orcid.org/0000-0002-6284-3033}{This research benefited from the support of the FMJH Program PGMO.}
\author{St.~Elmo Wilken}{Institute of Quantitative and Theoretical Biology, Heinrich Heine University Düsseldorf, Germany}{}{https://orcid.org/0000-0002-4113-2590}{This research benefited from the support of the FMJH Program PGMO.}
\author{Sebastian Pokutta}{Zuse Institute Berlin, Technische Universität Berlin, Germany}{}{https://orcid.org/0000-0001-7365-3000}{Research reported in this paper was partially supported through the Research Campus Modal funded by the German Federal Ministry of Education and Research (fund numbers 05M14ZAM,05M20ZBM) and the Deutsche Forschungsgemeinschaft (DFG) through the DFG Cluster of Excellence MATH+.}

\authorrunning{H.~Troppens et al.} 

\Copyright{Hannah Troppens, and Mathieu Besançon, and St.~Elmo Wilken, and Sebastian Pokutta} 



\category{} 

\relatedversion{} 

\supplement{Software (Source Code): \url{https://github.com/hannahtro/Loopless_Fluxes_with_Mixed_Integer_Optimization}}

\acknowledgements{}

\nolinenumbers 

\EventEditors{Petra Mutzel and Nicola Prezza}
\EventNoEds{2}
\EventLongTitle{23rd International Symposium on Experimental Algorithms (SEA 2025)}
\EventShortTitle{SEA 2025}
\EventAcronym{SEA}
\EventYear{2025}
\EventDate{July 22--24, 2025}
\EventLocation{Venice, Italy}
\EventLogo{}
\SeriesVolume{338}
\ArticleNo{7}

\usepackage[nocomma,short]{optidef}
\Crefname{constraint}{Constraint}{Constraints}
\crefname{constraint}{Constraint}{Constraints}
\usepackage{algorithm}
\usepackage{algpseudocode}
\usepackage[colorinlistoftodos,prependcaption,textsize=small]{todonotes}
\usepackage{makecell}

\begin{document}
\maketitle

\begin{abstract}
Constraint-based metabolic models can be used to investigate the intracellular physiology of microorganisms. 
These models couple genes to reactions, and typically seek to predict metabolite fluxes that optimize some biologically important metric.
Classical techniques, like Flux Balance Analysis (FBA), formulate the metabolism of a microbe as an optimization problem where growth rate is maximized. 
While FBA has found widespread use, it often leads to thermodynamically infeasible solutions that contain internal cycles (loops). 
To address this shortcoming, Loopless-Flux Balance Analysis (ll-FBA) seeks to predict flux distributions that do not contain these loops.
ll-FBA is a disjunctive program, usually reformulated as a mixed-integer program, and is challenging to solve for biological models that often contain thousands of reactions and metabolites. \\
In this paper, we compare various reformulations of ll-FBA and different solution approaches. 

Overall, the combinatorial Benders' decomposition is the most promising of the tested approaches with which we could solve most instances. However, the model size and numerical instability pose a challenge to the combinatorial Benders' method.
\end{abstract}

\section{Introduction}
In systems biology, detailed annotations of a microorganism's genome can be leveraged to construct a genome-scale metabolic model (GSMM). These models connect genes to proteins that catalyze specific reactions within the metabolism of a microbe. One essential and challenging task is predicting the \emph{flux} channeled by these metabolic reactions, i.e., the amount of molecules consumed by the different reactions. Reaction kinetics drive the flux of individual reactions, but due to challenges in their measurement and estimation, encoding them into models remains highly impractical, leaving them under-determined. Constraint-based models have been proposed to address these issues and improve flux estimation by reformulating the model as an optimization problem with reaction fluxes as variables.

In the following, we briefly introduce key biological concepts. The chemical reactions that take place in an organism, known as \textit{metabolism}, determine the central functions of the cell \cite{intro_computational_systems_biology}. Understanding metabolism is, therefore, crucial for understanding cellular behavior. \textit{Metabolites} are small molecules involved in biochemical reactions, and each \textit{reaction} converts a set of \textit{reactants} into \textit{products}, with defined \textit{stoichiometry} indicating the number of molecules involved. \textit{Enzymes} are important for cell metabolism as they \textit{catalyze} reactions, meaning they accelerate reactions without being consumed. We differentiate between \textit{internal reactions}, which involve only internal metabolites, and \textit{exchange reactions}, which interact with the environment. A \textit{metabolic network} captures this structure, often modeled as a directed hypergraph, where nodes represent metabolites and hyperedges represent reactions.

By assuming that the cell operates at metabolic steady state, the network of feasible reaction fluxes can be captured in a polytope.
Finally, the addition of a biologically relevant objective function, e.g. growth rate maximization, yields a linear program (LP) that can be solved efficiently.
This approach, coined Flux Balance Analysis (FBA), can be summarized in the following LP:
\begin{maxi}
  {\scriptstyle \mathbf v}{\mathbf c^\intercal \mathbf v}{\tag{FBA}\label{eq:FBA}}{}
    \addConstraint{\mathbf S \mathbf v= \mathbf 0}
    \addConstraint{\mathbf l \leq \mathbf v \leq \mathbf u.}
\end{maxi} 
The structure of the reaction network is captured in the stoichiometric matrix $\mathbf S \in \mathbb{R}^{m\times n}$, where $m$ is the number of internal metabolites and $n$ is the number of reactions.
The fluxes $\mathbf v \in \mathbb{R}^n$ are bounded by $\mathbf l, \mathbf u \in \mathbb{R}^n$. The linear system $\mathbf S \mathbf v= \mathbf 0$ ensures that the model is mass-balanced at steady state. 
FBA is inspired by viewing evolutionary selection pressure as a mathematical optimization problem to maximize fitness (biomass growth)~\cite{palsson_systems_biology}.
Importantly, FBA enabled systems biologists to predict cellular behavior in accordance with observations in experiments~\cite{FBA}.
However, solutions of~\ref{eq:FBA} often violate thermodynamic principles, as they may yield solutions that contain internal cycles, which are biologically unrealistic. An internal loop or internal cycle is a nonzero flux vector $\boldsymbol \ell$ such that the internal network is at steady-state: $\mathbf S_{\mathcal{I}} \boldsymbol \ell = \mathbf 0$ \cite{noor_proof_2012}. 
Loopless-Flux Balance Analysis (ll-FBA) is a constraint-based approach to predict flux distributions that removes internal loops~\cite{elimination_infeasible_loops}, thus optimizing over a subset of the FBA polytope.
Unlike FBA, ll-FBA is an $\mathcal{NP}$-hard disjunctive program and is challenging to solve in theory and practice for biological models that often contain thousands of reactions and metabolites.
In addition to tractability challenges, ll-FBA instances create numerical issues for solvers, motivating the design of robust and stable methods.

Besides ll-FBA, other methods exist to obtain loopless or thermodynamically feasible flux distributions. For example, thermodynamic Flux Balance Analysis incorporates thermodynamic constraints directly using metabolite concentrations, improving realism but requiring additional data \cite{tmfa}. Parsimonous FBA \cite{lewis_omic_2010} and CycleFreeFlux \cite{desouki_cyclefreeflux_2015-1} are post-processing approaches that remove cycles after solving FBA, offering computational efficiency but potentially sacrificing optimality.
 
In this paper, we tackle the ll-FBA problem with mixed-integer optimization techniques, aiming to bridge the gap between model biological realism and tractability even on genome-scale instances.
We derive disjunctive formulations of ll-FBA that notably avoid introducing some of the artificially large bounds on continuous variables.
We then design a combinatorial Benders' (CB) decomposition for the disjunctive problem, exploiting a natural separation between the flux and edge activation on one side, and thermodynamic feasibility on the other side.
We evaluate the different formulations and algorithms on realistic genome-scale metabolic model examples from the systems biology literature. To demonstrate the flexibility of our framework, we also conduct experiments incorporating enzyme constraints.
Computational experiments show that our CB algorithm performs much better than the alternatives and reveal choices of the algorithm's parameters that work well on the whole instance set, specifically studying the number of cuts to add per Benders' iteration and how to select these cuts. Our implementations of the various methods for ll-FBA are open-source and rely on open-source solvers and libraries to be directly utilized by the systems biology community.

\paragraph*{Notation}
A vector $\mathbf{v} = (v_1, ..., v_n) \in \mathbb{R}^n$ is identified as a column vector and printed in bold. 
With $\mathbf v \leq \mathbf w$ we denote element-wise inequality. 
The 1-vector is written as $\mathbf 1 := (1, 1, ..., 1) \in \mathbb{R}^n$, the 0-vector as $\mathbf 0 := (0, 0, ..., 0) \in \mathbb{R}^n$. The \textit{support} of $\mathbf{v}$ is the set of indices $i$ with $v_i \neq 0$ and is denoted by $\text{supp}(\mathbf v)$. With $\text{sign}(\mathbf v)$ the element-wise sign function is applied to $\mathbf v$.

The zero matrix is denoted by $\mathbf 0_{m,n}$ with $m$ rows and $n$ columns. With $\text{diag}(\mathbf v)$ we denote the quadratic matrix with $\mathbf v$ on the diagonal and 0 for all other entries.
A \textit{basis} $B$ of a vector space $V$ is a set of vectors $(\mathbf v_1, \mathbf v_2, ..., \mathbf v_n)$ that are linearly independent and every $\mathbf v \in V$ can be written as a linear combination of vectors in $B$.
The \textit{nullspace} of a matrix $\mathbf A \in \mathbb{R}^{m \times n}$ is defined as $\text{null}(\mathbf A):=\{\mathbf x \in \mathbb{R}^n: \mathbf A \mathbf x = ~\mathbf0\}$. 

\section{Loopless-FBA}
Loopless Flux Balance Analysis (ll-FBA) is a constraint-based approach that predicts flux distributions in cells that do not contain internal loops by incorporating thermodynamic information as an extension of FBA:
\begin{maxi}
    {\scriptstyle \mathbf v, \boldsymbol{\Delta^{\mu}}, \boldsymbol \mu}{\mathbf c^\intercal \mathbf v}{\label{eq:llfba_original}}{}
    \addConstraint{\mathbf S \mathbf v= \mathbf 0} 
    \addConstraint{\mathbf l \leq \mathbf v \leq \mathbf u}
    \addConstraint{\Delta^{\mu}_i v_i < 0 \lor v_i = 0  \quad \forall i \in \mathcal{I} \quad \quad \quad \quad} 
    \addConstraint{\boldsymbol{\Delta^{\mu}} = \mathbf S_{\mathcal{I}}^\intercal \boldsymbol \mu,} 
\end{maxi}
where $\mathbf v \in \mathbb{R}^n$, $\boldsymbol \mu \in \mathbb{R}^m$ and $\boldsymbol{\Delta^{\mu}} \in \mathbb{R}^{|\mathcal{I}|}$. The matrix $\mathbf S_\mathcal{I}$ is the submatrix of $\mathbf S$ that contains the columns of internal reactions only. $\boldsymbol \mu$ denotes Gibbs free energy, which predicts the direction of a reaction at constant pressure and temperature, and $\boldsymbol{\Delta^{\mu}}$ can be interpreted as the change in Gibbs free energy. A reaction $i$ is thermodynamically feasible only if the flux and the corresponding change in Gibbs free energy $\boldsymbol{\Delta^{\mu}_i}$ have opposite signs. Forward flux requires $\boldsymbol{\Delta^{\mu}_i} < 0$, while backward flux requires $\boldsymbol{\Delta^{\mu}} > 0$ \cite{beard_thermodynamic_2004}. However, in ll-FBA only the sign($\mathbf{\Delta^{\mu}}$) corresponds to the actual Gibbs free energy change. 
$\boldsymbol{\Delta^{\mu}} = \mathbf S_{\mathcal{I}}^\intercal \boldsymbol \mu$ ensures that the reaction energies around a circuit sum up to zero, which is known as \textit{Kirchhoff's second law}. 
Note that Kirchhoff's second law can also be ensured by the equation $\mathbf B^\intercal \boldsymbol{\Delta^{\mu}} = \mathbf 0$, where $\mathbf B \in \mathbb{R}^{|\mathcal{I}| \times k}$ is a matrix with columns formed by the set of vectors $\{\mathbf b_i\}_{i=1}^k$ that form a basis for $\mathrm{null}(\mathbf S_{\mathcal{I}})$ \cite{elimination_infeasible_loops}.
Since $\mathbf v$ and $\boldsymbol{\Delta^{\mu}}$ have different lengths, we assume that the first entries of $\mathbf v$ correspond to the internal reactions, so that each index $i$ refers to the pair $(v_i$, $\Delta^{\mu}_i)$.

Problem~\eqref{eq:llfba_original} is $\mathcal{NP}$-hard and much more complicated to solve than FBA~\cite{cornelis_metabolic_nodate}:
it is challenging due to the disjunction, the strict equality and the product of decision variables $\boldsymbol{\Delta^{\mu}}$ and $\mathbf v$.
To simplify the model, the disjunction is rewritten by introducing the constant $\epsilon$, and we arrive at the ll-FBA model we will use throughout the paper:
\begin{maxi}
    {\scriptstyle \mathbf v, \boldsymbol{\Delta^{\mu}}, \boldsymbol \mu}{\mathbf c^\intercal \mathbf v}{\tag{{ll-FBA}}\label{llFBA}}{}
    \addConstraint{\mathbf S \mathbf v= \mathbf 0} 
    \addConstraint{\mathbf l \leq \mathbf v \leq \mathbf u} 
    \addConstraint{\begin{aligned} &\bigl((v_i \geq 0) \land (\Delta^{\mu}_i \leq - \epsilon) \bigr) \, \, \lor \\ 
    &\bigl((v_i \leq 0) \land (\Delta^{\mu}_i \geq \epsilon) \bigr) \end{aligned} \quad \quad \forall i \in \mathcal{I}} 
    \addConstraint{\boldsymbol{\Delta^{\mu}} = \mathbf S_{\mathcal{I}}^\intercal \boldsymbol \mu.} 
\end{maxi} 
To avoid the degenerate solution where all $\Delta^{\mu}_i = 0$, it is standard practice to exclude the interval $[- \epsilon, \epsilon ]$ \cite{elimination_infeasible_loops}. Since only the sign of $\boldsymbol{\Delta^{\mu}}$ is thermodynamically relevant, this scaling does not alter the biological interpretation.
A solution to ~\ref{llFBA} does not contain internal cycles, but \ref{llFBA} is computationally harder as we need to solve a disjunctive program, which is usually reformulated as a mixed-integer program (MIP)\footnote{While ll-FBA is usually introduced in the literature directly via indicator constraints and their MIP reformulation using big-M constraints \cite{elimination_infeasible_loops}, we refer to this disjunctive form as ll-FBA throughout the paper and distinguish it from its reformulations used in practice.}.

\begin{figure}[h!]
    \centering
    \includegraphics[width=0.6\textwidth]{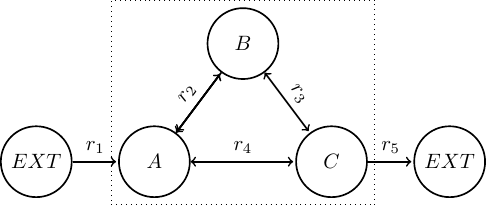}
    \caption{Simple model with internal loop}
    \label{fig:loop}
\end{figure}

As an example, let us consider the metabolic network in \cref{fig:loop}. We have three internal metabolites $A,B,C$, two irreversible exchange reactions $r_1, r_5$ and three reversible internal reactions $r_2, r_3, r_4$. The stoichiometric matrix is: 

\begin{equation*} \label{Eq:S_loop}
    \mathbf S =
    \left[\begin{array}{ccccc}
        1 & -1 & 0 & -1 & 0 \\
        0 & 1 & -1 & 0 & 0 \\
        0 & 0 & 1 & 1 & -1 \\
    \end{array}\right]
\end{equation*}

 and we assume the following bounds on the fluxes: 
\begin{equation*} \label{Eq:bounds_loop}
    \mathbf l = (0,-30,-30,-30,0) \quad \mathbf u = (10,30,30,30,10)
\end{equation*}

\begin{figure}[H]
    \centering
    \begin{subfigure}{0.49\textwidth}
        \includegraphics[width=\linewidth]{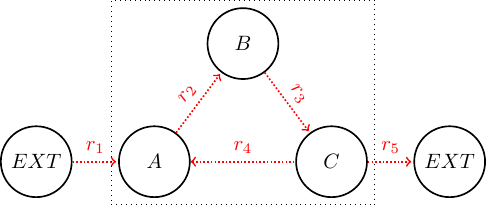}
        \caption{}
    \end{subfigure}%
    \hfill
    \begin{subfigure}{0.49\textwidth}
        \includegraphics[width=\linewidth]{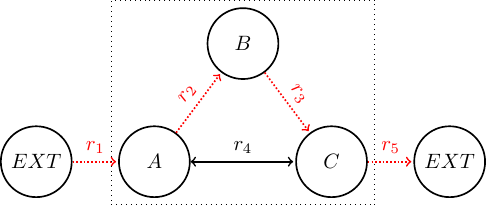}
        \caption{}
    \end{subfigure}
    \caption{Used reactions in (a) the \textsf{FBA} solution which contains an internal cycle, and (b) the \textsf{ll-FBA} solution.}
    \label{fig:loop_solutions}
\end{figure}

If we maximize the flux through the internal reactions, that is $\mathbf c = (0,1,1,1,0)$, an optimal solution is $\mathbf v^* = (10, 30, 30, -20, 10)$. The solution contains an internal loop, as there is a flux of 20 going through each of the internal reactions, as seen in \cref{fig:loop_solutions}.

\subsection*{Enzyme Information}
In classical FBA simulations, the predicted growth rate of the model is limited by the substrate uptake rate. 
This flux constraint on the uptake reaction is typically empirically measured, and imposed on the model. 
Recent theory has shown that the protein density in microbes is relatively constant, and thus a proteome limitation is a more plausible constraint to impose on the model~\cite{scott2010interdependence}.
Enzyme constrained metabolic models extend genome scale metabolic models by incorporating simplified enzyme kinetics.
In short, reaction flux is assumed to be proportional to a kinetic constant, and the concentration of enzyme available to catalyze a specific reaction. 
The optimization goal is then to distribute enzyme resources, subject to a total enzyme capacity bound, to maximize growth.
This extension can be used to model complex metabolic phenomena, like overflow metabolism.
Models incorporating this extra information are called enzyme constrained flux balance analysis models~\cite{improving_phenotype_predictions}.
 
Mathematically, if reaction $j$ is catalyzed by enzyme $E_i$, the reaction rate $v_j$ (in $\frac{\text{mmol}}{\text{gDW} \times \text h}$) depends on the enzyme usage $e_i$ (in $\frac{\text{mmol}}{\text{gDW}}$) and the turnover number (kinetic parameter) $k_{cat}^{ij}$ (in $\frac{1}{\text h}$) as $k_{cat}^{ij} e_i = v_j$, which is known as enzyme mass balance.
To account for enzyme mass balance in a GSMM, the stoichiometric matrix $\mathbf S$ is extended to form $\mathbf S^{ENZ}$. 
This extension adds a row for each enzyme $E_i$ and a column for the corresponding enzyme usage $e_i$, with $p$ enzymes. The lower left submatrix contains the enzyme information on the diagonal, that is $-1/k_{cat}^{ij}$. 
The resulting matrix is of the form:
\begin{equation*}
    \mathbf S^{ENZ} = \left[\begin{array}{ccc|ccc} 
        s_{1,1} & \text{...} & s_{1,n} & 0 & \text{...} & 0 \\ 
        \vdots & \ddots & \vdots & \vdots & \ddots & \vdots \\
        s_{m,1} & \text{...} & s_{m,n} & 0 & \text{...} & 0 \\
        \hline 
        -1/k_{cat}^{11} & \text{...} & 0 & 1 & \text{...} & 0 \\ 
        \vdots & \ddots & \vdots & \vdots & \ddots & \vdots \\
        0 & \text{...} & -1/k_{cat}^{pn} & 0 & \text{...} & 1
    \end{array}\right] = 
    \left[\begin{array}{cc} 
        \mathbf S & \mathbf 0_{m,p} \\
        \text{diag} (-1/k_{cat}^{ij}) & \mathbf I_p
    \end{array}\right].
\end{equation*}
In addition to the reaction rates $\mathbf v$, we are interested in the enzyme usage $\mathbf e$. We obtain the following optimization problem:
\begin{maxi}
    {\scriptstyle \mathbf v, \mathbf e}{\mathbf c^\intercal \mathbf v}{\tag{\texttt{ENZ}}\label{ENZ}}{}
    \addConstraint{\mathbf S^{ENZ} (\mathbf v, \mathbf e)= \mathbf 0} 
    \addConstraint{\mathbf l \leq \mathbf v \leq \mathbf u}
    \addConstraint{\mathbf 0 \leq \mathbf e \leq [\mathbf E],}
\end{maxi}
 where $[E_i]$ is the intracellular concentration of enzyme $E_i$.

 \section{Mathematical Models and Solving Methods}
In this section, we present various reformulations of the \ref{llFBA} problem in using indicator constraints and big-M constraints to model the disjunctions. Our goal is to compare the performance of 
\begin{enumerate}
\item solving these reformulations directly,  
\item decomposing the problem,
\item solving the convex hull formulation.
\end{enumerate}

\subsection{Big-M and Indicator Reformulations}
We can reformulate \ref{llFBA} to no longer write the disjunctions explicitly.
The mixed-integer programming formulation of flux balance analysis without unbounded internal cycles, using indicator constraints, is provided by \cite{elimination_infeasible_loops}:
\begin{maxi}
    {\scriptstyle \mathbf v, \boldsymbol{a}, \boldsymbol{\Delta^{\mu}}, \boldsymbol \mu}{\mathbf c^\intercal \mathbf v}{\tag{{ll-FBA (indicator)}}\label{ll-FBA (indicator)}}{}
    \addConstraint{\mathbf S \mathbf v=\mathbf 0} 
    \addConstraint{\mathbf l \leq \mathbf v \leq \mathbf u}    
    \addConstraint{a_i = 1}{\quad \implies \quad v_i \geq 0}{\quad \forall i \in \mathcal{I}} 
    \addConstraint{a_i = 1}{\quad \implies \quad \Delta^{\mu}_i \leq - \epsilon}{\quad \forall i \in \mathcal{I}}
    \addConstraint{a_i = 0}{\quad \implies \quad v_i \leq 0}{\quad \forall i \in \mathcal{I}} 
    \addConstraint{a_i = 0}{\quad \implies \quad \Delta^{\mu}_i \geq \epsilon}{\quad \forall i \in \mathcal{I}}
    \addConstraint{\boldsymbol{\Delta^{\mu}} = \mathbf S_{\mathcal{I}}^\intercal \boldsymbol \mu} 
    \addConstraint{a_i \in \{0,1\} \quad \forall i \in \mathcal{I},}
\end{maxi}
where each binary variable $a_i$ indicates which disjunct holds.
The mixed-integer program of flux balance analysis without unbounded internal cycles using big-M constraints, introducing the constant $M$, takes the form ~\cite{elimination_infeasible_loops}:
\begin{maxi}
    {\scriptstyle \mathbf v, \boldsymbol a, \boldsymbol{\Delta^{\mu}}, \boldsymbol \mu}{\mathbf c^\intercal \mathbf v}{\tag{{ll-FBA (big-M)}}\label{ll-FBA (big-M)}}{}
    \addConstraint{\mathbf S \mathbf v= \mathbf 0} 
    \addConstraint{\mathbf l \leq \mathbf v \leq \mathbf u}
    \addConstraint{-Ma_i + \epsilon(1-a_i) \leq \Delta^{\mu}_i \leq - \epsilon a_i + M(1-a_i) \quad \forall i \in \mathcal{I}}        
    \addConstraint{-M(1-a_i) \leq v_i \leq M a_i \quad \forall i \in \mathcal{I}}
    \addConstraint{\boldsymbol{\Delta^{\mu}} = \mathbf S_{\mathcal{I}}^\intercal \boldsymbol \mu}
    \addConstraint{a_i \in \{0,1\} \quad \forall i \in \mathcal{I}.} 
\end{maxi}
The big-M formulation ensures the opposite sign of $v_i$ and $\Delta^{\mu}_i$ if $v_i \neq 0$. If the flux through reaction $i$ is a forward flux, that is $v_i > 0$ and $a_i=1$, it holds that $-M \leq \Delta^{\mu}_i \leq -\epsilon$ and analogously for a backward flux. The value of $\epsilon$ affects the scale of $\boldsymbol{\Delta^{\mu}}$ (and $\boldsymbol \mu$). The value of $\epsilon$ should be bigger than the solver precision to differentiate it from 0. We set $\epsilon$ to 1 to match the formulation of \cite{elimination_infeasible_loops}. The big-M constant is the maximal absolute value of the flux bounds $\mathbf l, \mathbf u$.

\subsection{Decomposing ll-FBA}
Instead of solving a mixed-integer reformulation of \ref{llFBA} directly, the problem is decomposed into a master problem and a subproblem. The solving procedure involves first solving the master problem and then verifying the feasibility of the relaxed solution.
The master problem is a relaxed version of \ref{ll-FBA (indicator)}, where the loopless constraints are ignored, but the indicator variables $\boldsymbol a$ are assigned:
\begin{maxi}
    {\scriptstyle \mathbf v, \boldsymbol a}{\mathbf c^\intercal \mathbf v}{\tag{{MP}}\label{(MP)}}{}
    \addConstraint{\mathbf S \mathbf v= \mathbf 0} 
    \addConstraint{\mathbf l \leq \mathbf v \leq \mathbf u}
    \addConstraint{a_i = 1}{\quad \implies \quad v_i \geq 0}{\quad \forall i \in \mathcal{I}}      \addConstraint{a_i = 0}{\quad \implies \quad v_i \leq 0}{\quad \forall i \in \mathcal{I}.}
\end{maxi}
The indicator constraints can be linearized by using big-M constraints.
The subproblem is the following parameterised program:
\begin{maxi}
    {\scriptstyle \boldsymbol{\mu, \Delta^{\mu}}}{0}{\tag{{SP}}\label{(SP)}}{}
    \addConstraint{a_i^{MP} = 1}{\quad \implies \quad \Delta^{\mu}_i \leq - \epsilon}{\quad \forall i \in \mathcal{I}}
    \addConstraint{a_i^{MP} = 0}{\quad \implies \quad \Delta^{\mu}_i \geq \epsilon}{\quad \forall i \in \mathcal{I}}
    \addConstraint{\boldsymbol{\Delta^{\mu}} = \mathbf S_{\mathcal{I}}^\intercal \boldsymbol \mu,} 
\end{maxi}
where $\boldsymbol a^{\mathrm{MP}}$ are the values of a solution to the master problem.
If a relaxed solution is not a feasible solution to the \ref{llFBA} problem, a cut is added to the relaxed problem that cuts off the assignment of binary variables. This procedure is repeated until a relaxed solution is a feasible solution to \ref{llFBA}.

First, we introduce no-good cuts. If a relaxed solution $\boldsymbol a^{\mathrm{MP}}$ is thermodynamically infeasible, the following constraint is added to the master problem:

\begin{equation} \label{noGoodCut}
\sum_{i \in \mathcal{I}: \, a_i^{MP}=0} a_i + \sum_{i \in \mathcal{I}: \, a_i^{MP}=1} (1-a_i) \geq 1,
\end{equation}
where the entire infeasible relaxed solution is excluded from the solution space of the master problem.

\subsubsection*{Combinatorial Benders' Cuts}
When deriving combinatorial Benders' cuts for \ref{llFBA}, we aim to identify the variables that lead to thermodynamic infeasibility and potentially derive stronger cuts.
If the subproblem is infeasible, we compute a corresponding minimal infeasible subsystem (MIS) $\mathcal{C}$. The minimal infeasible subsystem contains a minimal number of constraint indices that make the subproblem infeasible.
The following combinatorial Benders' (CB) cut is added to the master problem if the subproblem is infeasible:
\begin{equation}
\sum_{i \in \mathcal{C}:\, a_i^{MP}=0} a_i + \sum_{i \in \mathcal{C}: \, a_i^{MP}=1} (1-a_i) \geq 1.
\end{equation}
Note that if $\mathcal{C}$ contains all internal reactions $\mathcal{I}$, the CB cut corresponds to a no-good cut (\cref{noGoodCut}).

\begin{algorithm}
    \caption{solving \ref{llFBA} with the combinatorial Benders' approach}\label{alg:CB}
    \begin{algorithmic}[1]
        \Require stoichiometric matrix $\mathbf S$, lower bound on fluxes $\mathbf l$, upper bound on fluxes $\mathbf u$, allowed number of cuts per iteration $m$, indices of internal reactions $\mathcal{I}$
        \State MP $\gets \textsf{build\_master\_problem}(\mathbf S, \mathbf l, \mathbf u, \mathcal{I})$
        \State $\mathbf{v, a} \gets \textsf{optimize(MP)}$ 
        \State $\mathcal{C} \gets \textsf{compute\_mis}(\mathbf S_\mathcal{I}, \boldsymbol a, m)$ \Comment{set of minimal infeasible subsystems}
        \While{$\mathcal{C} \neq \emptyset$}
            \State $\textsf{add\_cut}(\text{MP}, \boldsymbol a, \mathcal{C})$ \Comment{combinatorial Benders' cut}
            \State $\mathbf v, \boldsymbol a \gets \textsf{optimize}$(MP)
            \State $\mathcal{C} \gets \textsf{compute\_mis}(\mathbf S_\mathcal{I}, \boldsymbol a, m)$
        \EndWhile
        \State $\text{SP} \gets \textsf{build\_sub\_problem}(\mathbf S_\mathcal{I}, \mathcal{I}, \boldsymbol a)$
        \State $\boldsymbol{\Delta^{\mu}}, \boldsymbol \mu \gets \textsf{optimize}$(SP)
    \State \Return $(\mathbf v, \boldsymbol a, \boldsymbol{\Delta^{\mu}}, \boldsymbol \mu)$ \Comment{loopless flux distribution}
    \end{algorithmic}
\end{algorithm}

The pseudocode for the combinatorial Benders' cut procedure is shown in \cref{alg:CB}. 
The set of minimal infeasible subsystems is computed in \textsf{compute\_mis}, which identifies at most $m$ subsystems. 
The derivation of an MIS is based on the dual problem of the infeasible subproblem, and is explained in detail in \cite{codato_combinatorial_2006}. 

\noindent The primal problem in inequality form is:
\begin{maxi*}
    {\scriptstyle \mathbf x}{0}{\label[problem]{problem:MISPrimalStandard}}{} 
    \addConstraint{\mathbf{\tilde{A}} \boldsymbol \mu \leq \mathbf{\tilde{b}},}
\end{maxi*}
\quad where 
\begin{equation*}
    \mathbf{\tilde{A}} = 
    \left[\begin{aligned}
        & [\mathbf S_{\mathcal{I}}]_{*,i} \quad &\forall i \in \mathcal{I}: a_i^{MP} = 1\\
        - & [\mathbf S_{\mathcal{I}}]_{*,i} \quad &\forall i \in \mathcal{I}: a_i^{MP} = 0\\
    \end{aligned}\right]  
    \quad \mathbf{\tilde{b}} = \left[\begin{aligned}
        -&\epsilon^{|\mathcal{C}|} \\
    \end{aligned}\right].
\end{equation*} 

\noindent The linear program to find minimal infeasible subsystems is:
\begin{maxi}
    {\scriptstyle \boldsymbol \lambda}{\sum_i w_i \lambda_i}{\label{eq:mis_lp}}{} 
    \addConstraint{\mathbf{\tilde A}^\intercal \boldsymbol \lambda = \mathbf 0}
    \addConstraint{\boldsymbol \lambda \geq \mathbf 0}
    \addConstraint{\mathbf{\tilde b}^\intercal \boldsymbol \lambda = \mathbf 1,}
\end{maxi}
where $w_i$ is the weight corresponding to the dual variable $\lambda_i$. The support of each solution at a vertex of the feasible region of \eqref{eq:mis_lp}  defines a minimal infeasible subsystem. 
By modifying $\mathbf w$, we potentially derive several minimal infeasible subsystems to one infeasible solution $\boldsymbol a$.
Each solution corresponds to a minimal infeasible subsystem. The nonzero elements in $\boldsymbol \lambda$
correspond to a MIS with the set of reaction indices $\mathcal{C}$. We choose how many cuts $k$ can be added per iteration relative to the model size. For any $i \in {1, ..., k}$, we set the i-th coefficient in the objective function to zero and set all other coefficients to 1.
At the end of the MIS search, we filter the minimal infeasible subsystems found to have unique subsystems within each iteration.
To account for the model size, we set $k$ to a percentage of the number of reactions of the model. For instance, MIS $0.5\%$ means that at most $0.5\%$ minimal infeasible subsystems are blocked per iteration.
Instead of using all minimal infeasible subsystems found, we can filter the subsystems and block only a subset in the master problem (see \cref{section:appendix_cut_selection}).

\subsection{Convex Hull Reformulation}
The feasible region of the convex hull reformulation corresponds to the convex hull of the disjunctive constraints and thus represents the tightest possible convex relaxation of the original disjunctive program \cite{balas1979disjunctive}. However, this formulation introduces a significantly larger number of decision variables and constraints compared to the big-M approach. In particular, the number of additional constraints grows exponentially with the number of disjunctions.
Each disjunction is associated with two binary variables:
\begin{align}
    \left[\begin{array}{cc} y_i \\ v_i \geq 0 \\ \Delta^{\mu}_i \leq -\epsilon \end{array} \right] \lor 
        \left[\begin{array}{cc} y_{i + |\mathcal{I}|} \\ v_i \leq 0 \\ \Delta^{\mu}_i \geq \epsilon \end{array} \right] 
        \quad \forall i \in \mathcal{I}.
\end{align}
To formulate the disjunctive constraint with a mixed-integer program,
the constraint below is added to ensure the disjunction holds:
\begin{equation*}
    y_i + y_{i + |\mathcal{I}|} = 1 \quad \forall i \in \mathcal{I}.
\end{equation*}
Two continuous variables are added for each decision variable in a disjunction. The variables corresponding to $v_i$ are denoted by $v_{i_1}$ and $v_{i_2}$ and the variables corresponding to $\Delta^{\mu}_i$ are $\Delta^{\mu}_{i_1}$ and $\Delta^{\mu}_{i_2}$. The hull reformulation of the disjunction is:
\begin{align*}
    v_i &= v_{i_1} + v_{i_2} \\
    \Delta^{\mu}_i &= \Delta^{\mu}_{i_1} + \Delta^{\mu}_{i_2} \\
    - v_{i_1} &\leq 0 \\
    v_{i_2} &\leq 0 \\
    \Delta^{\mu}_{i_1} &\leq - y_i \\ 
    - \Delta^{\mu}_{i_2} &\leq -y_{i + |\mathcal{I}|}. 
\end{align*}

\section{Computational Experiments}
In this section, we present the results of the computational experiments for the \ref{llFBA} problem. We compare solving \ref{ll-FBA (big-M)} and \ref{ll-FBA (indicator)} with \texttt{SCIP} directly, to solving the convex hull reformulation and to decomposing the problem with a combinatorial Benders' approach.

\subsection{Experimental Setup}
We implemented the different solving methods\footnote{available at \\{\fontsize{9}{48} \url{https://github.com/hannahtro/Loopless_Fluxes_with_Mixed_Integer_Optimization}}} in Julia 1.9.0. We use \texttt{JuMP} \cite{JuMP} and \texttt{MathOptInterface} \cite{mathoptinterface} to build the mathematical models. The \texttt{DisjunctiveProgramming} package is used to get the hull formulation. The MIP solver used in the experiments is \texttt{SCIP} 9.0 ~\cite{scip8,scip9}.
The LP solver used is \texttt{HiGHS} 1.8.1 \cite{HiGHS}. We use \texttt{COBREXA} \cite{cobrexa} to load the biological model data.
The experiments were carried out on a 64-core compute node equipped with an Intel Xeon Gold 6338 2GHz CPU and 512GB RAM. The CPU memory was limited to 3000MB, and the time limit is 1800 seconds unless specified otherwise.

\subsubsection*{Metabolic Model Instances}
We test the entire set of metabolic networks of the \textit{biochemical, genetic, and genomic}~(BiGG) database\footnote{last accessed on Sep 25, 2024} \cite{BiGG}. The models cover different model sizes, ranging from the smallest model e\_coli\_core with 72 metabolites and 95 reactions to the largest model Recon3D with 5835 metabolites and 10600 reactions.

We use a subset of metabolic networks of yeast models made available by \cite{lu_yeast_2021}. The models selected are similar in size to the largest BiGG models used here (see \cref{Tab:yeast_model_size}).

When adding enzyme constraints to the model, the flux rate depends on the enzymes and the flux bounds $\mathbf l, \mathbf u$ are only required to restrict the direction of a reaction. Due to the limited availability of curated enzyme data in the BiGG database, we use randomly generated enzyme parameters for the BiGG models. Our aim is not to ensure biological realism in these cases, but rather to demonstrate the flexibility of our framework in handling enzyme constraints.
We use \texttt{COBREXA} to build the enzyme models, and the enzyme data (i.e. the turnover numbers and protein molar masses) is generated randomly as specified below. The turnover numbers differentiate for the forward and backward direction of a reaction and therefore we split each reversible reaction into one forward and one backward reaction.
At least one enzyme is mapped to each reaction including the turnover numbers for the reaction in the forward and the turnover number for the backward direction are defined. 
The turnover numbers are taken independently from a standard normal distribution. 
The protein concentrations have to be in the interval $[0, 1000]$. 
An enzyme is made up of one or more proteins.
Each protein is randomly associated with either mass group A or B, and the product mass of each group is bounded by 0.5 $\frac{\mathrm{mmol}}{\mathrm{gDW}}$.
For each protein we assign a product molar mass randomly from a uniform distribution.  

\subsection{Results}
\begin{figure}[h!]
\centering
\includegraphics[width=1.0\textwidth]{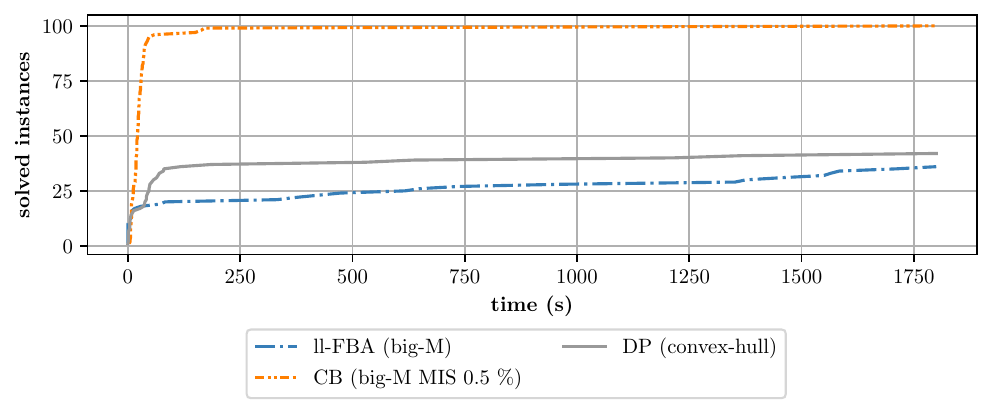}
\caption{\label{fig:dp_performance} Comparing the number of BiGG instances optimally solved by the big-M and convex hull formulations, by solving \ref{ll-FBA (big-M)} directly, and by solving the problem with CB (CB (big-M MIS 0.5 \%)). Note that some instances remain unsolved by the combinatorial Benders' approach, as shown in \cref{Tab:termination_CB} and \cref{Tab:termination_mis}.}
\end{figure}

We begin by comparing the performance of three approaches for solving \ref{ll-FBA (big-M)} for the BiGG instances: solving it directly with \texttt{SCIP}, solving its convex hull reformulation, and applying combinatorial Benders' cuts.
\cref{fig:dp_performance} shows that the convex hull reformulation solves the fewest instances.
Even though the feasible region of the hull reformulation is the convex hull, the reformulation introduces more constraints and decision variables, resulting in increased running times.  
In contrast, when comparing the MIP reformulations to solving \ref{ll-FBA (big-M)} with combinatorial Benders' cuts, we observe that the latter achieves superior performance. Specifically, it solves the majority of instances and does so significantly faster.

We next present detailed computational results for decomposing \ref{llFBA} with a particular focus on combinatorial Benders's cuts.
We compare the performance of solving \ref{ll-FBA (big-M)} and \ref{ll-FBA (indicator)} with \texttt{SCIP} directly on the BiGG instances to two decomposition approaches: solving the decomposition of \ref{ll-FBA (big-M)} with no-good cuts and solving it with the combinatorial Benders' cuts, all within a time limit of 1800 seconds. For the combinatorial Benders' method, we experiment with three different formulations of the master problem. In each iteration of \cref{alg:CB}, we add a single cut. The formulations we consider are: 
\begin{enumerate}
    \item the indicator formulation,
    \item the big-M formulation,
    \item a combined approach using both the indicator and big-M formulations (i.e.~adding both groups of constraints).
\end{enumerate} 

\begin{figure}[h!]
\centering
\includegraphics[width=0.7\textwidth]{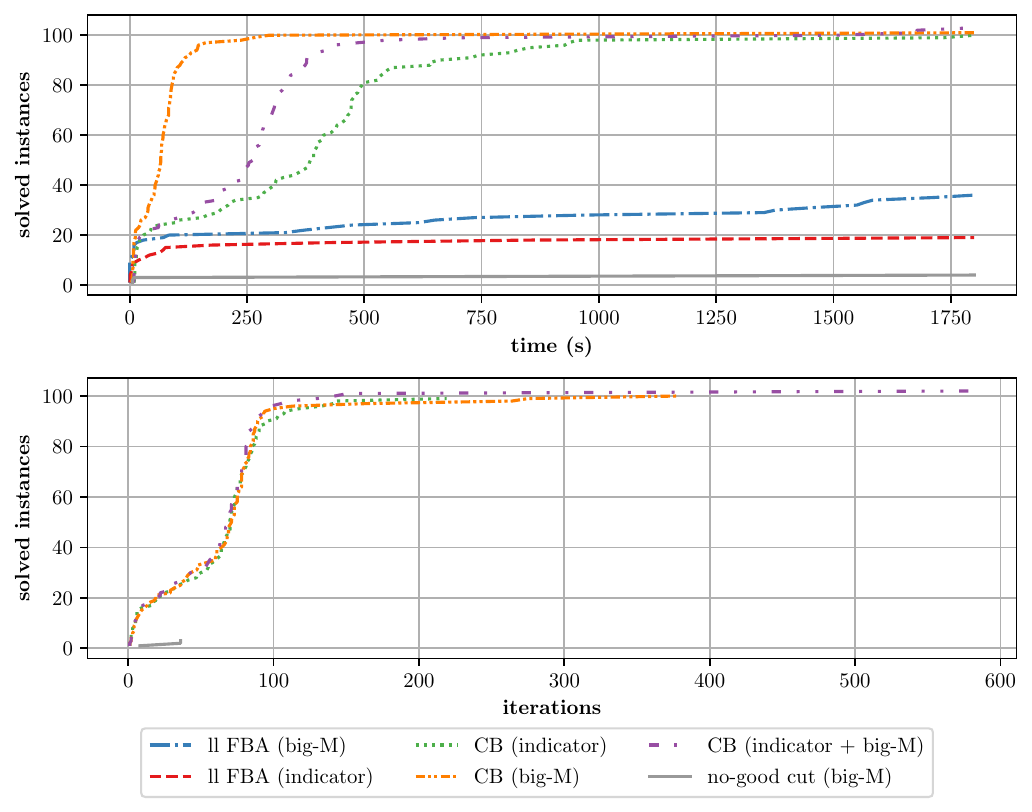}
\caption{Performance of CB master problem variants to directly solving ll-FBA}
\label{fig:comparison_solved_instances_indicator_and_big_m_as_MP}
\footnotesize{Comparing the number of optimally solved BiGG instances of \ref{ll-FBA (big-M)}, \ref{ll-FBA (indicator)} directly, solving the decomposition with no-good cuts and solving the combinatorial Benders' decomposition of ll-FBA within the time limit of 1800 seconds. We experiment with different formulations of the master problem: big-M reformulation, indicator constraints, and indicator and big-M constraints. Per iteration of the combinatorial Benders' approach, one minimal infeasible subsystem is computed.}
\end{figure}

\cref{fig:comparison_solved_instances_indicator_and_big_m_as_MP} shows that the combinatorial Benders' method outperforms the other approaches, solving most instances within the time limit. The fastest combinatorial Benders' method is the big-M formulation, followed by the combined indicator and big-M formulation, and last the indicator formulation.
Although the running time varies among the different formulations, the number of iterations remains similar across methods, suggesting that the running time is primarily determined by the time spent on solving the master problem. The combinatorial Benders' method with indicator constraints solves the fewest instances among the tested setups, yet we continue our experiments with it, because it solves a different set of instances to optimality compared to the method with big-M constraints.

In contrast, directly solving any of the ll-FBA variants results in significantly fewer optimally solved instances within the time limit. The poorest performance is observed with the no-good cuts approach, which solves only $8\%$ of the instances.
We see that combinatorial Benders' cuts are much stronger than no-good cuts.
This poor performance was expected since no-good cuts only remove exactly one solution per iteration, in contrast to Benders' cuts which block a cycle and thus produces a constraint that removes a large portion of the binary hypercube.

\cref{Tab:termination_CB} shows that the combinatorial Benders' method with big-M constraints solves most instances (93 \%) while needing the shortest average solving time (58 seconds).  
For the combinatorial Benders' decomposition of \mbox{\ref{ll-FBA (big-M)}}, the instances \texttt{iAM\_Pc455}, \texttt{iAM\_Pf480} and \texttt{iCN718} error during the solution process, even though solving \mbox{\ref{ll-FBA (big-M)}} returns an optimal solution. If we solve \texttt{iAM\_Pc455} and \texttt{iAM\_Pf480} with a tighter tolerance, we find an optimal solution within the time limit. The combinatorial Benders' method with indicator constraints also solves 93\% of the instances to optimality, but the average solving time is 227 seconds. 

On the other hand, directly solving any of the ll-FBA variants leads to significantly fewer optimally solved instances within the time limit. Specifically solving \ref{ll-FBA (indicator)} directly leads to 17\% optimally solved instances, while solving \ref{ll-FBA (big-M)} directly leads to 33\% optimally solved instances. 
Instances \texttt{iSB619} and \texttt{iCN718} are said to be infeasible when solving \ref{ll-FBA (indicator)} directly and are therefore excluded from further experiments with indicator constraints.

However, in instances where solving \ref{ll-FBA (big-M)} takes less than 10 seconds, it achieves an average solving time of just one second, outperforming all other methods. This shows there is a small overhead to designing dedicated algorithms compared to passing a single model to a MIP solver.

\begin{table}[!ht]
    \centering
    \begin{tabular}{@{\extracolsep{4pt}}lrrrrrrrrr@{}}
    \hline
        \multicolumn{2}{c}{} & \multicolumn{2}{c}{\textbf{\thead{ll-FBA \\ (big-M)}}} & \multicolumn{2}{c}{\textbf{\thead{ll-FBA \\(indicator)}}} & \multicolumn{2}{c}{\textbf{\thead{CB \\ (big-M)}}}  & \multicolumn{2}{c}{\textbf{\thead{CB \\ (indicator)}}} 
        \\ \cline{1-4} \cline{5-6} \cline{7-8} \cline{9-10}
        \thead{solved \\in (s)} & \thead{\# \\ inst.} & \thead{\% \\ solved} & \thead{time\\ (s)} & \thead{\% \\ solved} & \thead{time\\ (s)} & \thead{\% \\ solved} & \thead{time\\ (s)} & \thead{\% \\ solved} & \thead{time\\ (s)} \\ \hline
        0-10 & 22 & 91 & \textbf{8} & 64 & 72 & \textbf{95} & 10 & \textbf{95} & 19 \\
        10-60 & 23 & 43 & 1139 & 4 & 1442 & \textbf{100} & \textbf{42} & \textbf{100} & 201 \\
        60-100 & 44 & 7 & 1671 & 7 & 1531 & \textbf{100} & \textbf{80} & 95 & 505 \\
        100-600 & 12 & 17 & 1542 & 0 & 1800 & \textbf{100} & \textbf{137} & 92 & 639 \\
        600-1800 & 6 & 0 & 1800 & 0 & 1800 & 0 & 1800 & \textbf{33} & \textbf{996} \\ \hline
        all & 107 & 33 & 517 & 17 & 829 & \textbf{93} & \textbf{58} & \textbf{93} & 227 \\ \hline
    \end{tabular}
    \caption{\label{Tab:termination_CB} \small Comparing the number of optimally solved BiGG instances and the average running time of solving \ref{ll-FBA (indicator)} directly, of solving \ref{ll-FBA (big-M)} directly and solving their combinatorial Benders' decomposition within a time limit of 1800 seconds. Per iteration of the combinatorial Benders' approach, one minimal infeasible subsystem is computed. The BiGG instances are divided depending on the shortest running time required to solve \ref{ll-FBA (big-M)} with any of the listed methods. The time column shows the average time in seconds using the geometric mean.}
\end{table} 

\begin{figure}[htb]
\centering
\begin{subfigure}{0.49\textwidth}
\centering
\includegraphics[width=1\textwidth]{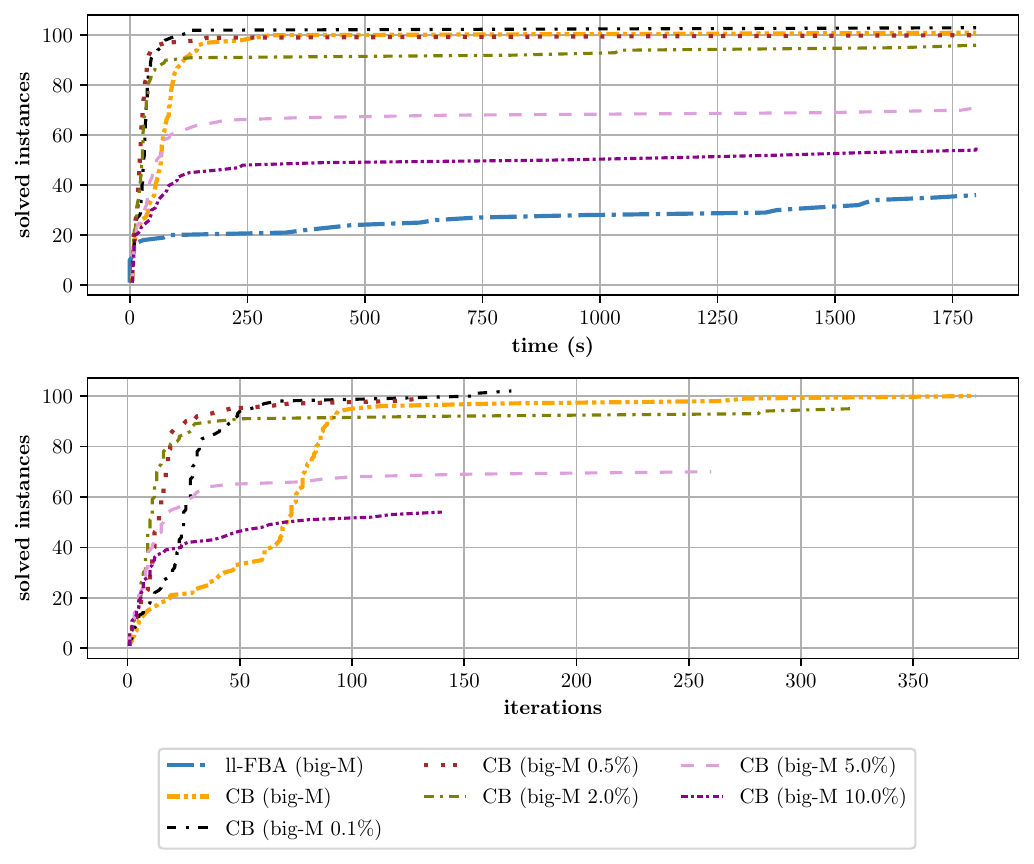}
\caption{CB with big-M.}
\label{fig:cbbigmmultiple}
\end{subfigure}
\begin{subfigure}{0.49\textwidth}
\centering
\includegraphics[width=1\textwidth]{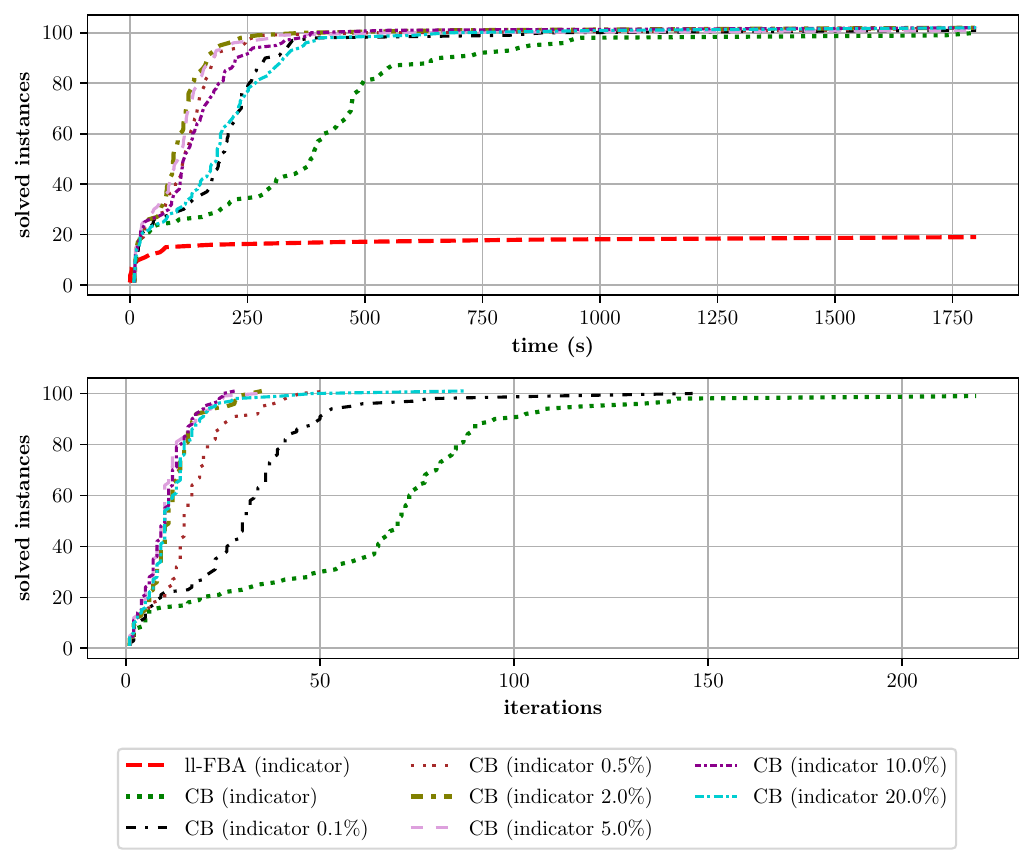}
\caption{CB with indicator.}
\label{fig:decomp1}
\end{subfigure}
\caption{Performance of CB with multiple cuts, comparing the number of optimally solved BiGG instances of solving ll-FBA (indicator/big-M) directly and running CB.
We experiment with the number of cuts added per iteration of CB depending on the instance size.
MIS 0.5\% means that at most $k$ cuts are added per iteration, where $k$ is 0.5\% of the number of reactions of the model.}
\label{fig:multicuts}
\end{figure}

So far, we added one combinatorial Benders' cut per iteration of \cref{alg:CB}.
In \cref{fig:multicuts}, we compare the performance of directly solving ll-FBA to the combinatorial Benders' method for both indicator and big-M constraints, experimenting with the number of cuts added per iteration.
In the following, all potential cuts found per iteration are added. We investigate different cut selection strategies in \ref{section:appendix_cut_selection}.

As expected, increasing the number of cuts per iteration leads to reduced running times. Any combinatorial Benders' setup using multiple cuts per iteration outperforms the single-cut approach. 
However, setups using big-M constraints and adding MIS 2\% or more cuts per iteration perform significantly worse than combinatorial Benders' with one cut per iteration, with fewer instances being solved within the time limit. 
On the contrary, when using indicator constraints, the CB algorithm maintains its capacity to solve most instances across all choices of the number of cuts per round. 

\cref{Tab:termination_mis} highlights that the combinatorial Benders' decomposition \mbox{\ref{ll-FBA (big-M)}} achieves the best performance, solving most instances in the shortest average time. Specifically, blocking up to 0.1\% of the number of reactions in the per iteration allows to solve 95\% of the instances, while needing only around $28$ seconds in average for the solving process when blocking up to 0.5\% of the number of reactions.
However, adding too many cuts per iteration can negatively impact the overall performance, leading to more iterations and longer running times. For instance, blocking up to 5\% of the number of reactions per iteration results in an average running time of $130 \text{ seconds}$, which is higher than the performance of solving it with indicator constraints $(86 \text{ seconds})$ and higher than adding just one cut per iteration $(58 \text{ seconds})$. Additionally, the instance \texttt{iAM\_Pv561} is said to be infeasible, even though we find a solution with other combinatorial Benders' setups. If we tighten the tolerance of the solving process, the instance is solved to optimality showing that too many cuts can increase numerical instability.  

\begin{table}[!ht]
    \centering
    \begin{tabular}{@{\extracolsep{4pt}}rrrrr@{}}
    \hline
    \multicolumn{1}{c}{} & \multicolumn{2}{c}{\textbf{CB (indicator)}} & \multicolumn{2}{c}{\textbf{CB (big-M)}} \\
    \cline{2-3} \cline{4-5}
    \thead{MIS \%} & \thead{\% solved} & time (s) & \thead{\% solved} & time (s) \\ \hline
    0.1 & 93 & 137.74 & \textbf{95} & \textbf{30.91} \\ 
    0.5 & \textbf{94} & 96.11 & 93 & \textbf{27.55} \\ 
    2.0 & \textbf{94} & 84.06 & 90 & \textbf{38.63} \\ 
    5.0 & \textbf{94} & \textbf{86.0} & 67 & 129.19 \\ 
    10.0 & \textbf{94} & \textbf{101.57} & 50 & 275.82 \\ \hline
    CB & \textbf{93} & 227.12 & \textbf{93} & \textbf{57.75} \\ \hline
    \end{tabular}
    \caption{Comparing the number of optimally solved instances and average running time of solving \ref{ll-FBA (indicator)} and \ref{ll-FBA (big-M)} with the combinatorial Benders' method on the BiGG instances within a time limit of 1800 seconds. We experiment with the number of cuts added per iteration of the combinatorial Benders' approach depending on the instance size. MIS 0.5\% means that at most $k$ cuts are added per iteration, where $k$ is 0.5\% of the number of reactions of the model. The time column shows the geometric mean. The last row shows the results of the CB methods adding one cut per iteration.}
    \label{Tab:termination_mis}
\end{table}

After experimenting with various formulations of the master problem and different numbers of cuts per iteration in \cref{alg:CB}, we test both solving strategies on yeast models, which are larger than most models of the BiGG database. Neither solving strategy performs well on the yeast instances: only one strategy is able to solve a single instance out of the 11 instances tested. 

\begin{figure}[h!]
    \centering
    \includegraphics[width=0.7\textwidth]{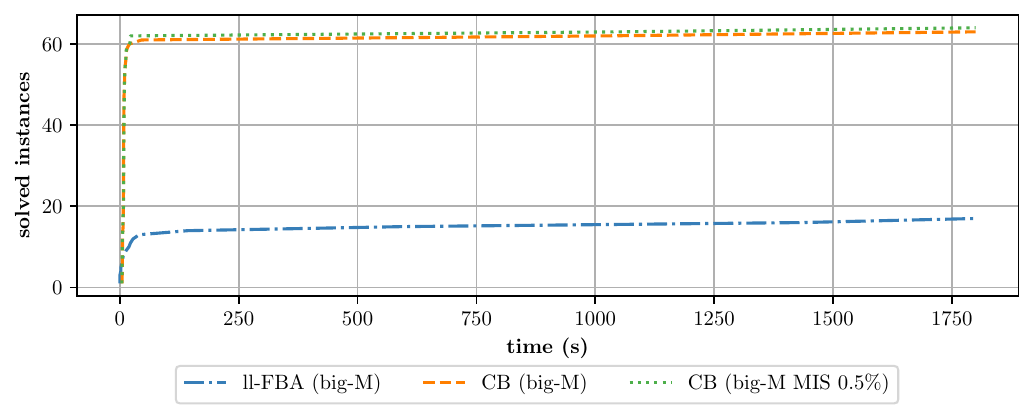}
    \caption{\label{fig:comparison_gecko} Performance of the combinatorial Benders' method on enzyme models.
    Comparing the number of optimally solved instances with enzyme constraints of solving \ref{ll-FBA (big-M)} directly, and solving the problem with the combinatorial Benders' approach within the time limit of 1800 seconds. We experiment with adding a single cut per iteration of the combinatorial Benders' approach and with adding at most $k$ cuts per iteration, where $k$ is 0.5\% of the number of reactions in the model.}
\end{figure}

Finally, we experiment with solving enzyme models using the combinatorial Benders' approach. We compare directly solving \ref{ll-FBA (big-M)} to solving it with the combinatorial Benders' decomposition, as the big-M formulation showed the best performance on previous experiments. We compare blocking one cycle per iteration to blocking $k$ cycles per iteration, where $k$ is 0.5\% of the total number of reactions in the model. The results are presented in \cref{fig:comparison_gecko}.
The results indicate that both variants of the combinatorial Benders' approach perform similarly and outperform solving \ref{ll-FBA (big-M)} directly. Both combinatorial Benders' variants require only a few iterations for most instances. As the reaction rate is limited by enzyme abundance, the relaxed solution already contains fewer loops, requiring fewer cuts to obtain an optimal solution.

However, the combinatorial Benders' approach is less stable for enzyme models. While around 70\% of instances are solved to optimality, 18\% of instances have their solution process terminated due to numerical errors.
Indeed, the enzyme models are numerically challenging, most errors arise from the master problem at the first iteration not resulting in the same objective value as \ref{eq:FBA} within a tolerance of $10^{-3}$.  
In our experiments, we are able to solve most enzyme instances within a few seconds.
Since the enzyme data is randomly generated, we do not make any claims on the general performance on enzyme models.
However, the results show the flexibility of the combinatorial Benders' approach, which can easily be extended with additional constraints.
They also hint that enzyme constraints regularize FBA models and remove some thermodynamically infeasible solutions from the original polytope.

\section{Conclusion}

In this paper, we compared the performance of different formulations of \ref{llFBA}.
We compared solving different reformulations of \ref{llFBA} with the MIP solver \texttt{SCIP} directly including the big-M reformulation and the convex hull formulation. 
We then experimented with different solution methods in order to solve more instances than directly solving ll-FBA.
We experimented with decomposing the ll-FBA problem into a master and a subproblem, and compared separation approaches using no-good cuts and combinatorial Benders' cuts.
The Combinatorial Benders' approach performs much better on BiGG instances than the big-M formulation of ll-FBA.
However, the performance of the combinatorial Benders' approach depends strongly on the number of cuts added per iteration: too many cuts per iteration slow down the solution process significantly.
On the yeast models, which are larger than most BiGG models, neither with the big-M reformulation nor with most combinatorial Benders' setups are we able to solve instances to optimality, showing the need for further work on scalable algorithms for ll-FBA.
The combinatorial Benders' approach is flexible and can be extended, as long as we can still decompose the problem into a master and a subproblem, we illustrate this with additional enzymatic constraints.
The models with enzyme constraints are experimentally much faster to solve. This opens the interesting question of their tractability: under suitable assumptions, is the ll-FBA with enzyme constraints hard? 



\bibliography{bibliography}

\newpage
\appendix

\section{Model Data}
\begin{table}[!ht]
    \centering
    \addtolength{\leftskip} {-2cm}
    \addtolength{\rightskip}{-2cm}
    \small
    \begin{tabular}{lll}
    \hline
        \textbf{model} & \textbf{\# metabolites} & \textbf{\# reactions} \\ \hline
        \texttt{Hanseniaspora\_uvarum} & 2464 & 3569 \\
        \texttt{yHMPu5000035696\_Hanseniaspora\_singularis} & 2460 & 3534 \\
        \texttt{yHMPu5000034963\_Hanseniaspora\_clermontiae} & 2464 & 3573 \\
        \texttt{yHMPu5000035695\_Hanseniaspora\_pseudoguilliermondii} & 2476 & 3559 \\
        \texttt{yHMPu5000035684\_Kloeckera\_hatyaiensis} & 2465 & 3582 \\
        \texttt{Eremothecium\_sinecaudum} & 2549 & 3471 \\
        \texttt{yHMPu5000035659\_Saturnispora\_dispora} & 2598 & 3378 \\
        \texttt{Tortispora\_caseinolytica} & 2693 & 3597 \\
        \texttt{Starmerella\_bombicola\_JCM9596} & 2695 & 3735 \\
        \texttt{Eremothecium\_gossypii} & 2556 & 3555 \\
        \texttt{Ashbya\_aceri} & 2553 & 3623 \\ \hline
    \end{tabular}
    \caption{\label{Tab:yeast_model_size} Model size of used yeast models.}
\end{table}

\section{Additional Computational Experiments}
\subsection{Experiments with Yeast Instances}
\begin{table}[!ht]
\centering
\begin{tabular}{@{\extracolsep{4pt}}lrrr@{}}
    \hline
    \textbf{\thead{solving strategy}} & \textbf{\thead{\# optimal}} & \textbf{\thead{\# time limit}} & \textbf{\thead{\# error}} \\ \hline
    ll-FBA (big-M) & 0 & 11 & 0 \\
    CB (big-M) & 0 & 11 & 0 \\
    CB (indicator) & 0 & 11 & 0 \\
    CB (big-M MIS 0.1 \%) & 0 & 11 & 0 \\
    CB (big-M MIS 0.2 \%) & 0 & 11 & 0 \\
    CB (big-M MIS 0.5 \%) & 1 & 10 & 0 \\
    CB (big-M MIS 2.0 \%) & 0 & 11 & 0 \\
    CB (indicator MIS 0.1 \%) & 0 & 11 & 0 \\
    CB (indicator MIS 0.2 \%) & 0 & 11 & 0 \\
    CB (indicator MIS 0.5 \%) & 0 & 11 & 0 \\
    CB (indicator MIS 2.0 \%) & 0 & 10 & 1 \\ \hline
\end{tabular}
\caption{\label{Tab:yeast} Comparing the termination status of solving yeast instances with \ref{ll-FBA (big-M)}  directly and solving \ref{ll-FBA (big-M)} with CB.
We experiment with the indicator and big-M formulations of the master problem. We compare adding a single cut per iteration of CB and adding multiple cuts. MIS 0.5\% means that at most $k$ cuts are added per iteration, where $k$ is 0.5\% of the number of reactions of the model.}
\end{table}

\subsection{Cut Selection}\label{section:appendix_cut_selection}
Instead of blocking the entire set of cycles found per iteration of te combinatorial Benders' approach, we experiment with selecting a subset of cycles. We experiment with different strategies: 

\begin{enumerate}
    \item blocking distinct cycles,
    \item selecting the $k$ smallest cycles, 
    \item selecting cycles with a density limit.
\end{enumerate}

\begin{figure}[htb]
\centering
\begin{subfigure}{0.49\textwidth}
\centering
\includegraphics[width=1\textwidth]{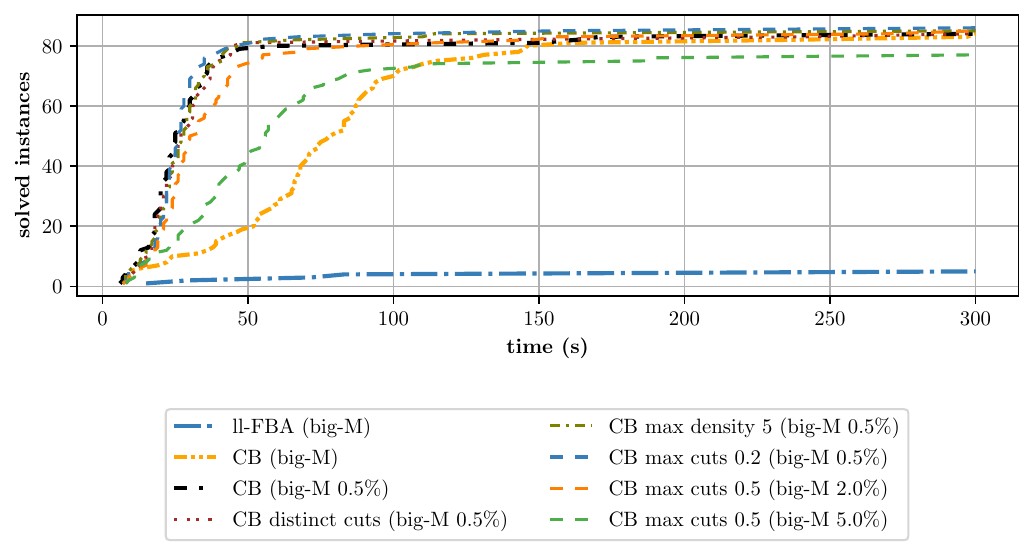}
\caption{CB with big-M.}
\label{fig:cbbigmselection}
\end{subfigure}
\begin{subfigure}{0.49\textwidth}
\centering
\includegraphics[width=1\textwidth]{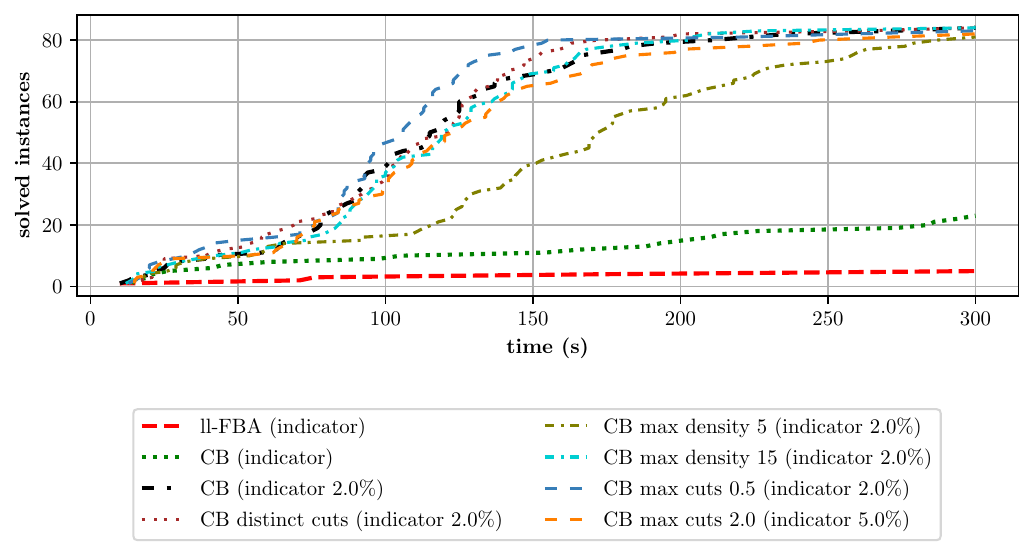}
\caption{CB with indicator.}
\label{fig:cbindicatorselection}
\end{subfigure}
\caption{Performance of CB (big-M/indicator) with different cut selection strategies.
We experiment with adding distinct cuts, adding the cuts with a maximal density and adding at most $k$ cuts (max cuts).
All strategies are applied with different parameters and a different size of maximally detected number of MIS per iteration.
}
\label{fig:cut_selection_instances}
\end{figure}

\cref{fig:cut_selection_instances} shows that all tested strategies outperform blocking one cycle per iteration and directly solving \ref{llFBA}.
Most strategies perform similarly regarding the running time and the number of optimally solved instances.
However, deriving an excessive number of MIS per iteration leads to increased running time and fewer solved instances.
In the case of CB with indicator constraints, since more time is spent solving the master problem compared to CB with the big-M formulation, the time spent on cut selection impacts the overall running time.
We observe that being too aggressive in cut filtering, i.e.~removing useful cuts, can increase the running time.

\end{document}